\documentclass[12pt]{article}
\usepackage{epsfig}
\usepackage{amsmath}

\usepackage{amssymb}
\newtheorem{Proposition}{Proposition}
  \newtheorem{Remark}{Remark}

  \newtheorem{Theorem}[Proposition]{Theorem}

\newenvironment{proof}{\par\noindent {\sc Proof.}}{$\Box$}

\newcommand{\be}{\begin{equation}}
\newcommand{\ee}{\end{equation}}

\def\z{\noindent}

\def\Box{{\hfill\hbox{\enspace${\sqre}$}} \smallskip}
\def\sqr#1#2{{\vcenter{\vbox{\hrule height .#2pt
                             \hbox{\vrule width .#2pt height#1pt \kern#1pt
                                   \vrule width .#2pt}
                             \hrule height .#2pt}}}}
\def\sqre{\mathchoice\sqr54\sqr54\sqr{4.1}3\sqr{3.5}3}

\def\CC{\mathbb{C}}
\def\RR{\mathbb{R}}
\def\ZZ{\mathbb{Z}}
\def\NN{\mathbb{N}}

\begin{document}

\title{Analyzability in the context of PDEs and applications}
\author{{\it O. Costin}\\ \small{Math Department, Rutgers University}\\
{\it S. Tanveer} \\
\small{Math Department, The Ohio State University}}

\maketitle
\date{}
\bigskip

\begin{abstract}
 \z  We discuss the notions of resurgence, formalizability, and formation of
  singularities in the context of partial differential equations. The results
  show that \'Ecalle's how analyzability theory extends naturally to PDEs.
\end{abstract}

\section{Introduction} 

The study of nonlinear partial differential equations in the complex domain
and especially formation of spontaneous singularities of their solutions is
not a well understood subject.  The theory of \'Ecalle's analyzable functions,
originally developed (mainly) for functions of one variable, provides a set of
tools which are well suited to address some of these issues, but the extension
to several variables is not immediate.

In the case
of linear ODEs under suitable assumptions, there is a complete system of
formal solutions as transseries, \cite{Ecalle}, and these are generalized
Borel (multi)summable to a fundamental system of actual solutions of the
system; sufficiently powerful results of a similar nature have been shown in
the nonlinear case as well \cite{Ecalle}, \cite{{BBRS}},  \cite{Invent},
\cite{DMJ}.  This is true for difference equations as well,
\cite{ Kuik, Braaksma, Braaksma-Kuik}.

While transseries solutions can be considered and their Borel summation shown
in the context of PDEs there are a number of difficulties specific to several
variables. We first discuss a number of specific obstacles to extending the
theory in a straightforward way, and then refer to what we expect to be a
general approach to many problems and overview a number of recent results
utilizing this approach.

\section{Difficulties of formalizability and analyzability in PDEs}

\subsection{Insufficiency of formal representations} For PDEs even the notion a general formal solution appears to elude a
definition that is reasonably simple and useful.

\z {\em Example 1.} The equation $f_t+f_x=0$ has the general solution
$$f(t,x)=F(t-x)$$
with $F$ any differentiable function; it is not clear to us
that a worthwhile definition can be associated to the description "general
formal expression in $t-x$"; on the other hand, restricting ourselves to
special, well defined, combinations in $t-x$ would correspondingly limit the
number of associated actual functions, precluding a complete solution of the
original PDE. In this example, actual solutions outnumber by far formalizable
ones.

\bigskip

\z The (apparently) opposite situation is possible as well.

\subsection{Absence of general summation procedures}

\z Let now $f_0\in C^\infty(a,b)$. The initial value problem

\begin{equation}
  \label{CR1}
  f_t+if_x=0;\ \ f_{\displaystyle
  |_{t=0}}=f_0(x)
\end{equation}

\z always has formal solution as $t\downarrow 0$:

\begin{equation}
  \label{fsCR}
\tilde{f}=\sum_{k=0}^\infty t^k
\frac{(-i)^kf_0^{(k)}(x)}{k!}
\end{equation}
but it has no actual solution, if $f_0$ is real-valued non-analytic (cf. the
proof of next proposition).
 
There is {\bf no nontrivial} {\em summation procedure of formal Taylor series
  over an interval} nor a more restricted one that would associate actual
solutions to (\ref{CR1}) as the following proposition shows. (See also Remark
\ref{R1}).

\begin{Proposition}\label{nosum}
  Let $S$ be a summation procedure defined on a differential algebra $D_S$ of
  formal series of the form

$$\tilde{g}=\sum_{k=0}^{\infty} g_k(x) t^k$$

\z where $S$ is assumed to have the following (natural) properties:
\begin{itemize}
\item $S$ is linear.

\item $S$ commutes with differentiation.
  
\item $S(\tilde{g})\sim g_0(x)$ as $t\downarrow 0$.
  
\item If $\tilde{f}\in D_S$ then $S\tilde{f}:(a,b)\times
  (0,\epsilon)\mapsto\CC$ (where $\epsilon$ is allowed to depend on $\tilde{f}$).
\end{itemize}

\z Assume if $f_0$ is real valued. Then $\tilde{f}$ in (\ref{fsCR}) is in
$D_S$ {\bf iff} $\tilde{f}$ is convergent (in the usual sense).
\end{Proposition}

\begin{proof}
  Assume $\tilde{f}\in D_s$. Then, by the properties of $S$, the function
  $f=S\tilde{f}$ is a differentiable function, and it is a solution of
  (\ref{CR1}) in a domain $\mathcal{D}=(a,b)\times (0,\epsilon)$.  If we write
  $f=u+iv$ we see that the pair $(u,v)$ satisfies the Cauchy-Riemann equations
  in $\mathcal{D}$ and thus $f$ is analytic in $z=t+ix$ with
  $(x,t)\in\mathcal{D}$. The third property of $S$ shows that $f_0$ is the
  limit as $z$ approaches the interval $(a,b)$ from the upper-half plane.
  Since $f_0$ is real valued, then by the Schwarz reflection principle $f$
  extends analytically through $(a,b)$ to a neighborhood of $(a,b)$ in the
  lower half plane; in particular $f$ is analytic on $(a,b)$. But then
  (\ref{CR1}), which is the Taylor series of $f$ at points on $(a,b)$ is
  convergent. \end{proof}

\subsection{Obstacles to determining the formal solutions} We now contrast 
formal analysis of ODEs and PDEs. \begin{enumerate}\item Consider the
  Painlev\'e P1 equation
  \begin{equation}
    \label{P1}
   y''=6y^2+x
  \end{equation}
  a rather nontrivial example of a second order nonlinear ODE. A detailed
  analysis of transseries and their Borel summability in this example are
  discussed in \cite{Invent}. We only mention a few aspects relevant to the
  present discussion.
  
  {\em Finding formal solutions of (\ref{P1}) is quite straightforward}.
  Searching first for algebraic behavior, dominant balance shows that
  $6y^2\sim -x$, say $y\sim 6^{-\frac{1}{2}}i {x}^{\frac{1}{2}}$, and then,
  consistent with this, $y''=o(x)$. It follows that a formal series expansion
  can be gotten by taking $y_0=0$ and then, for $n\in\NN,$ iterating the
  recurrence

$$y_{n+1}=6^{-\frac{1}{2}}i\sqrt{x-y_n''}$$

\z A power series solution is readily obtained in this way,

\begin{equation}
  \label{ps}
  \tilde{y}=6^{-\frac{1}{2}}i\sqrt{x}-\frac{1}{48x^2}+\frac{49\sqrt{6}
  i}{4608x^{9/2}}+...
\end{equation}

\z which is not classically convergent but is Borel summable to an actual
solution \cite{DMJ}; the complete transseries can be calculated and Borel
summed in a similar way, \cite{DMJ}. The possibility (and convenience) of the
formal calculation is partly due to {\em asymptotic simplification}, resulting
in a dominant balance equation, $$6y^2+x=0$$
which can be solved exactly, from
which a complete solution of the full problem follows by appropriate
perturbation theory.

\item 
Compare this problem to the periodically forced Schr\"odinger equation

\begin{equation}
  \label{eq:eqa} i\, \frac{\partial\psi }{\partial
t}\, =\left[-\partial_{xx}+V(x)+\Omega(x)\, \cos \omega t\right]\, \psi
\end{equation}
\z Under physically reasonable assumptions $\psi$ is transseriable \cite{UAB}:

\begin{Theorem}[\cite{UAB}]\label{genth1}
  
  Assume $\Omega, V$ are compactly supported and continuous, and $\Omega>0$
  throughout the support of $V$. For $t>0$ there exist $N\in\NN$ and
  $\{\Gamma_k\}_{k\le N}$, $\{F_{\omega;k}(t,x)\}_{k\le N}$,
  $2\pi/\omega$-periodic functions of $t$, such that

\begin{equation}\label{transP}\psi(t,x)=\sum_{j\in\ZZ}e^{ij\omega
t}h_j(t,x)+\sum_{k=1}^N e^{-\Gamma_k t}F_{\omega;k}(t,x)\end{equation}

\z with $\Re\Gamma_k>0$ for all $k\le N$, and
$h_j(t,x)=O\left(t^{-\frac{3}{2}},|j|!^{-\frac{1}{2}}\right)$ have Borel
summable power series in $t$,

\begin{equation}\label{e112}h_j(t,x)=\mathcal{LB}\sum_{k\ge 3}h_{kj}(x)t^{-k/2}
\end{equation}
The operator $\mathcal{LB}$ (Laplace-Borel) stands for generalized
Borel summation \cite{DMJ}.
\end{Theorem}

\z Insofar as a formal analysis would be concerned, it is to be noted that
there is no small parameter in (\ref{eq:eqa}) and largeness of $t$ does not
make any term negligibly small; a posteriori, knowledge of the transseries
(\ref{transP}) confirms this. In this sense, (\ref{eq:eqa}) admits no further
simplification.  There is, to the knowledge of the authors, no straightforward
formal way based on (\ref{eq:eqa}) to determine whether $\psi\to 0$, let alone
its asymptotic expansion.

\end{enumerate}
\section{Overcoming these difficulties: the approach of asymptotic
  regularization}

First note an implication of \'Ecalle's analyzability techniques
\cite{Ecalle-book}: a wide class of problems can be regularized by suitable
Borel transforms.  Summability of general solutions of ODEs or difference
equations, \cite{Ecalle-book, BBRS, IMRN, DMJ} shows that, under appropriate
transformations, the resulting equations admit {\em convergent} solutions, an
indication of the regularity of the associated equation.

Transseries are obtained, by suitable inverse transforms, from these
regularized solutions. 

{\em In the case of PDEs it appears that regularizing the equation is in many
  cases the adequate approach.} The result (\ref{transP}) is obtained in this
way.

\subsubsection{Elementary illustration: regularizing the heat equation}
\begin{equation}
  \label{heat}
   f_{xx}-f_t=0
\end{equation}
  
\z Since (\ref{heat}) is
parabolic,  power series solutions 
\begin{equation}
  \label{fsol}
  f=\sum_{k=0}^{\infty} t^k F_k(x)=\sum_{k=0}^{\infty} \frac{F_0^{(2k)}}{k!}t^k
\end{equation}

\z are divergent even if $F_0$ is analytic (but not entire).  Nevertheless,
under suitable assumptions, Borel summability results of such formal solutions
have been shown by Lutz, Miyake, and Sch\"afke \cite{LMS} and more general
results of multisummability of linear PDEs have been obtained by Balser
\cite{Balser}.

$\bullet$ The heat equation can be regularized by a suitable Borel summation.
The divergence implied, under analyticity assumptions, by (\ref{fsol}) is
$F_k=O(k!)$ which indicates Borel summation with respect to $t^{-1}$. Indeed,
the substitution

\begin{equation}
  \label{subst}
  t=1/\tau;\ \  f(t,x)=t^{-1/2}g(\tau,x)
\end{equation}
yields

$$ g_{xx}+\tau^2 g_{x}+\frac{1}{2}\tau g =0$$

\z which becomes after formal inverse Laplace transform (Borel transform) in
$\tau$,

\begin{equation}
  \label{wave0}
  p{\hat{g}}_{pp}+\frac{3}{2}{\hat{g}}_p-{\hat{g}}_{xx}=0
\end{equation}

\z which is brought, by the substitution
${\hat{g}}(p,x)=p^{-\frac{1}{2}}u(x,2p^{\frac{1}{2}});\ 
y=2p^{\frac{1}{2}}$, to the wave equation, which is hyperbolic, thus
{\em regular}
\begin{equation}
  \label{wave1}
u_{xx}-u_{yy}=0.
\end{equation}
Existence and uniqueness of solutions to regular equations is
guaranteed by Cauchy-Kowalevsky theory. For this simple equation the
general solution is certainly available in explicit form:
$u=u_-(x-y)+u_+(x+y)$ with $u_-,u_+$ arbitrary $C^2$
functions. Since the solution of (\ref{wave1}) is related to a
solution of (\ref{heat}) through (\ref{subst}), to ensure that we do
get a solution it is easy to check that we need to choose $u_-=u_+=u_0$ (up to
an irrelevant additive constant which can be absorbed into $u_-$) which
yields,

\begin{equation}
  \label{hk1}
 f(t,x)=t^{-\frac{1}{2}}\int_0^{\infty}y^{-\frac{1}{2}}\left[u_0 \left( x+2\,y^{\frac{1}{2}}
  \right)+u_0 \left( x-2\,y^{\frac{1}{2}} \right)\right]e^{ -y/t}dy
\end{equation}

\z which, after splitting the integral and making the substitutions 
$x\pm 2\,y^{\frac{1}{2}}=s$ is transformed into the usual heat kernel solution,

\begin{equation}
  \label{hk2}
  f(t,x)= t^{-\frac{1}{2}}\int_{-\infty}^{\infty}u_0(s)\exp\left(-\frac{(x-s)^2}{4t}\right)ds
\end{equation}
{\bf In conclusion} although there is perhaps no systematic way to
formalize the general solution of the heat equation, appropriate
inverse Laplace transforms allow us a complete solution of the problem
(in an appropriate class of initial conditions which ensure
convergence of the integrals).
\begin{Remark}\label{R1}
  Proposition \ref{nosum} can be also understood in the following way.
  Equation (\ref{CR1}) is already regular. {\em Any} actual solution, if it
  exists with the initial condition given in the Proposition, is trivially
  formalizable since it is then analytic. It is thus natural that no further
  summable formal solutions exist.
\end{Remark}

\subsection{Nonlinear equations: regularization by Inverse
  Laplace Transform} In this section we briefly mention a number of our
results of that substantiate regularizability.
\subsubsection{} Consider the third order scalar evolution PDE:
\begin{equation}
  \label{origeq}
f_t - f_{yyy} =\sum_{j=0}^3 b_j(y,t;f) f^{(j)} \,+\,r(y,t);~f(y, 0) =f_I
\end{equation}
Formal Inverse Laplace Transform with respect to $y$
gives\footnote{For technical convenience, in \cite{CPAM} we used
  oversummation. The paper \cite{Inventiones2} shows that in fact
  Borel summability holds in the correct variable, in the more general
  setting decribed in the next section.}
\begin{equation}
  \label{regeq}
  F_t +p^3 F =\sum_{j\le 3;k<\infty}\left [
B_{j,k} * (p^j F) * F^{*k} \right ] +R(p,t)
\end{equation}

\noindent where convolution is the Laplace one, $(f*g)(p)=\int_0^p
f(s)g(p-s)ds$.  This equation is regular in that formal power series in $p$
converges, since the coefficients in the equation are analytic.

Multiplying by the integrating factor of the l.h.s. and
integrating yields
\begin{multline*}
  F(p,t) = {\mathcal{N}} F (p,t)\\= F_0(p, t) +
\sum_{j\le 3;k<\infty} \int_{0}^t (-1)^j e^{-p^3
(t-\tau)} \left [ (p^j F)*B_{j,k}*F^{*k} \right ] (p, \tau) d\tau
\end{multline*}

\z The regularity of this equation plays a crucial role in the proofs in
\cite{CPAM} where we find the actual solutions of equation (\ref{origeq}).
\subsubsection{} Similar methods were later  extended  \cite{Inventiones2} to
equations of the form
$$
  {\bf u}_t + \mathcal{P}(\partial_{\bf x}^{\bf j}){\bf u}+{\bf g}
  \left ( {\bf x}, t, \{\partial_{\bf x}^{{\bf j}} {\bf u}\} \right )
  =0;\ {\bf {u}}({\bf x}, 0) ={\bf {u}}_I({\bf x})$$
  with
  $\mathbf{u}\in\CC^{r}$, for $ t\in (0,T)$ and large $|{\bf x}|$ in a
  poly-sector $S$ in $\mathbb{C}^d$ ($\partial_{\bf x}^{\bf j} \equiv
  \partial_{x_1}^{j_1} \partial_{x_2}^{j_2} ...\partial_{x_d}^{j_d}$
  and $j_1+...+j_d\le n$). The principal part of the constant
  coefficient $n$-th order differential operator $\mathcal{P}$ is
  subject to a cone condition. The nonlinearity ${\bf g}$ and the
  functions $\mathbf u_I$ and $\mathbf u$ satisfy analyticity and decay
  assumptions in $S$.

   The paper \cite{Inventiones2}  shows existence and uniqueness of the solution of this
  problem and finds its asymptotic behavior for large $|\bf x|$.

Under further regularity conditions on $\mathbf g$ and $\mathbf u_I$ which
ensure the existence of a formal asymptotic series solution for large
$|\mathbf x|$ to the problem, we prove its Borel summability to the actual
solution $\mathbf u$.
  
\smallskip

    In special cases motivated by applications we show how the method
  can be adapted to obtain short-time existence, uniqueness and
  asymptotic behavior for small $t$, of sectorially analytic solutions,
  without size restriction on the space variable.
\subsection{ Nonlinear Stokes
  phenomena and movable singularities} \z In the context of ODEs it was shown
\cite{Invent}, under fairly general assumptions, that the information
contained in the regularized problem (equivalently, in the transseries) can be
used to determine more global behavior of solutions of nonlinear equations, in
particular the fact that they form spontaneous singularity close to anti-Stokes
lines. The method, transasymptotic matching, was extended to difference
equations \cite{Braaksma, Braaksma-Kuik, Kuik}.

In nonlinear partial differential equations, formation of singularities is a
very important phenomenon but no general methods to address this issue
existed.

The method of regularization that we described provides such a method. We
briefly discuss the main points of \cite{CTprep}.

At present our methods apply to nonlinear evolution PDEs with {\em one space
  variable}; even for these, substantial new difficulties arise with respect
to \cite{Invent}.

Consider the modified Harry Dym equation (arising in Hele-Shaw dynamics)
\begin{eqnarray*}
\label{0.2.1}
H_t-H^3H_{xxx}+H_x-\frac{1}{2}H^3=0; \ \ H(x,0) = \frac{1}{\sqrt{x}}
\end{eqnarray*}
\z in an appropriate sector.

\z {\bf Small time behavior}. From \cite{CPAM} it follows that there exists a
unique solution to above problem, and it has Borel summable series for small
$t$ and small $y=x-t$:
 \begin{multline*}
   H(x,t) = y^{-1/2} - t \left ( \frac{15}{8 y^5} +\frac{1}{ 2 y^{3/2}}\right
   ) + t^2 \left (\frac{25875}{128 y^{19/2}}+ \frac{195}{32 y^6} + \frac{3}{8
       y^{5/2}}\right ) ~+~...
\end{multline*}

\z \z {}{{\bf Singularity manifolds near anti-Stokes lines}}.  To apply the
method of transasymptotic matching, we look on a scale where the asymptotic
expansion becomes formally invalid: $y=x-t = O(t^{2/9})$.  The transition
variable is thus
\begin{equation*}
\label{0.2.4}
\eta =\frac{x-t}{t^{2/9}}, ~~\tau = t^{7/9}, ~~H(x,t) = t^{-1/9}~ G (\eta, \tau )
\end{equation*}
Substituting into (\ref{0.2.1}), we obtain the following equivalent equation

\begin{equation*}
\label{0.2.5}
-\frac{G}{9} - \frac{2}{9} \eta G_\eta + \frac{7}{9} \tau G_\tau + 
\frac{\tau}{2} G^3 - G^3 G_{\eta \eta \eta} =0
\end{equation*}
\z The natural formal expansion solution in this regime is

\begin{equation}
\label{0.2.6}
G(\eta, \tau) = \sum_{k=0}^\infty \tau^k G_k (\eta)
\end{equation} 
\z with matching conditions at large $\eta$, to ensure the solution agrees
with the one obtained in \cite{CPAM}:
\begin{equation*}
\label{0.2.11}
G_0 (\eta) ~\sim ~\eta^{-1/2};\ \ G_k (\eta) ~\sim ~\frac{A_k}{\eta^{k+1/2}}
\end{equation*}

\z We show that the series (\ref{0.2.6}) is actually {\em convergent} and
equals $H(x,t)$ in the Borel summed region (the radius of convergence shrinks
however with $\eta$).  The convergence problem is subtle and required a rather
delicate construction of suitable invariant domains.  Having shown that, it is
intuitively clear (and not difficult to prove) that if $G_0$ is singular, then
$H$ is singular.  The leading order solution $G_0$ satisfies
\begin{equation*}
\label{0.2.7}
 G_0 + 2\, \eta\, G_0^\prime  +9\, G_0^3\, G_0^{\prime \prime \prime} = 0 
\end{equation*}
while for $k \ge 1$, 
\begin{equation*}
\label{0.2.8} 
G_0^3 {\mathcal L}_k G_k = R_k 
\end{equation*}
where
\begin{equation*} 
\label{0.2.9}
\mathcal{L}_k u = u^{\prime\prime\prime} + \frac{2}{9 G_0^3} \eta u^\prime - \left (\frac{7k-1}{9 G_0^3}
+ \frac{3 G_0^{\prime \prime \prime}}{G_0} \right ) u  
\end{equation*}
and the right hand side $R_k $  is given by
\begin{eqnarray*}
\label{0.2.10}
R_k (\eta) = \frac{1}{2} \sum_{k_1+k_2+k_3 = k-1} G_{k_1} G_{k_2}
G_{k_3} + \sum_{k_j < k, \sum k_j = k} G_{k_1} G_{k_2} G_{k_3} G_{k_4}^{\prime \prime \prime} 
\end{eqnarray*}
\z The nonlinear ODE of $G_0$ with asymptotic condition has been studied in
\cite{93} and computational evidence suggested clusters of
singularities $\eta_s$, where $G_0 (\eta)\sim e^{\pi i/3}\left (\frac{\eta_s}{3}
\right )^{1/3} (\eta - \eta_s)^{2/3} $.

\begin{equation}
  \label{ode1}
  G_0 + 2\, \eta\, G_0^\prime +9\, G_0^3\, G_0^{\prime \prime
    \prime} = 0 
\end{equation}
For a rigorous singularity analysis of (\ref{ode1}) we now used
transasymptotic matching as developed for ODEs \cite{Invent}. The asymptotic
behavior of $G_0$ is of the form
\begin{equation*}
\label{0.3.1}
G_0 (\eta) ~\sim~ \eta^{-1/2} 
U (\zeta )   
\end{equation*}
where
{\begin{equation*}
\label{0.3.1.5}
\zeta = -\ln C +\frac{9}{8}\ln \eta +  \frac{i 4\sqrt{2}}{27} \eta^{9/4} +(2n-1)\pi i
\end{equation*}}
(where $C$ is the Stokes constant) and $U(\zeta) $ satisfies algebraic
equation {\begin{equation*}
\label{0.3.2}
\frac{1}{4} e^{\zeta+2} = e^{-2 \sqrt{U}} ~\left ( \frac{\sqrt{U} +1}{\sqrt{U}-1} \right )
\end{equation*}}
Singularities of $ U(\zeta)$ occur at 
$\zeta_s =-\ln 4 + 2 -i \pi $, corresponding for $n\in\NN$ large
to
{\begin{equation*}
\label{0.3.2.1}
\frac{i 4\sqrt{2}}{27} \eta_s^{9/4} +\frac{9}{8} \ln ~\eta_s = -2 + \ln 4 - (2 n -1) i \pi+\ln C
\end{equation*}}
The {\textbf{Theorem} that we prove in \cite{CTprep} is that {\em
    For a singularity $\hat{\eta}_s$ of $G_0$ \footnote {with $|\hat{\eta}_s|$ large
      enough and with $\arg ~\hat{\eta}_s$ close to the anti-Stokes line $\arg ~\eta
      = - \frac{4 \pi}{9}$.} there exists a domain $\mathcal{D}$ around the
    singularity \footnote{that extends to $\infty$ with $\arg ~\eta \in \left
        ( -\frac{2\pi}{9}+\delta , \frac{2\pi}{9} -\delta \right )$ for some
      $\frac{\pi}{9} ~>~\delta ~>~0$, and includes a region $\mathcal{S}$
      around the the singularity $\hat{\eta}_s$ but excludes an open neighborhood.}
    such that the expansion is convergent for small $\tau$.}}

{\em In particular, for small $\tau$, the singularity of $G(\eta,
  \tau)=H(x,t)$ approaches the singularity of $G_0 $ and is, to the
  leading order, of the same type, $(\eta - \eta_s)^{2/3}$.}

\end{document}